# Stability of Collinear Equilibrium Points in Robe's Generalised Restricted Three Body Problem

**K.T. Singh$^a$, B.S. Kushvah$^b$ and B. Ishwar$^c$**

(*a*) Lecturer in Mathematics M.B. College, Imphal (Manipur),
(*b*) J.R.F D.S.T. Project*,
(*c*) Principal Investigator D.S.T. Project*

*University Department of Mathematics,
B.R.A. Bihar University Muzaffarpur-842001,India
Email: ishwar_bhola@hotmail.com

**Abstract:** We have examined the stability of collinear equilibrium points in Robe's generalized restricted three body problem. The problem is generalised in the sense that more massive primary has been taken as an oblate spheroid. We have found the position of equilibrium points. We have obtained variational equations of the problem. With the help of characteristic roots, we conclude that collinear points are unstable. Robe's result may be verified from the generalised result.

## 1. Introduction

"Robe's Restricted Three Body Problem", a model considered by Robe (1977), is a new kind of Restricted Three Body Problem in which one of the primaries is rigid spherical shell $m_1$ filled with a homogeneous incompressible fluid of density $\rho_1$. The second primary is a mass point $m_2$ outside the shell and the third body $m_3$ is a small solid sphere of density $\rho_3$ inside the shell with the assumption that the mass and radius of $m_3$ are infinitesimally small. He has shown the existence of an equilibrium point with $m_3$ at the centre of the shell while $m_2$ describes a Keplerian orbit around it. Further, he has discussed the linear stability of the equilibrium point for the whole range of parameters occurring in the problem.

Later on, several authors have studied the Robe's problem. Shrivastava and Garain (1991) have studied the effect of small perturbations in the Coriolis and centrifugal forces on the location of equilibrium points in the Robe's problem. Plastino and Plastino (1995) have considered the Robe's problem by taking the shape of the fluid body as Roches' ellipsoid. They have studied the linear stability of the equilibrium solution too. Giordano, Plastino and Plastino (1997) have discussed the effect of drag force on the existence and stability of the equilibrium points in Robe's problem. In all the above problems, the existence of only one equilibrium point is discussed. Hallan and Rana (2001) established that in the Robe's elliptic restricted three body problem, there is only one equilibrium point while in the Robe's circular restricted three body problem, there are two, three and infinite number of equilibrium points depending upon the certain values of parameters occurring in the problem.

In this paper we have studied the stability of the libration point in Robe's restricted three body problem. The problem is generalised in the sense that bigger primary is taken as an oblate spheroid.

## 2. Stability of Collinear Equilibrium Points

The equations of motion in he Robe's generalised restricted three body problem are given as

$$\ddot{x} - 2n\dot{y} = \frac{\partial \Omega}{\partial x}, \qquad \ddot{y} + 2n\dot{x} = \frac{\partial \Omega}{\partial y}, \qquad \ddot{z} = \frac{\partial \Omega}{\partial z} \qquad \text{---- (1)}$$

where,

$$\Omega = \frac{n^2}{2}\left(x^2 + y^2\right) - kr_1 + \frac{\mu}{r_2} \qquad \text{---- (2)}$$

$$r_1^2 = \left(x + \mu\right)^2 + y^2 + z^2, \qquad r_2^2 = \left(x + \mu - 1\right)^2 + y^2 + z^2 \qquad \text{---- (3)}$$

The angular velocity is given by

$$n^2 = 1 + \frac{3}{2}A_1, \quad k = \frac{4}{3}\pi\,\rho_1\left(1 - \frac{\rho_1}{\rho_3}\right), \quad \mu = \frac{m_2}{m_1 + m_2}$$

Coordinates of $m_1$ and $m_2$ are $(-\mu, 0, 0)$ and $(1-\mu, 0, 0)$. The libration points exist when $\Omega_x = \Omega_y = \Omega_z = 0$

We have found the collinear libration point at

$$x = \left(1 - \frac{3}{2}A_1 + 4k\right)\mu, \qquad y = 0, \qquad z = 0 \qquad \text{---- (4)}$$

At equilibrium points

$$\Omega_{xy}^0 = \Omega_{yx}^0 = \Omega_{zx}^0 = \Omega_{xz}^0 = 0 \qquad \text{---- (5)}$$

and $\qquad \Omega_{zz}^0 = -2k - \frac{\mu}{r_2^3} \qquad \text{---- (6)}$

$$\Omega_{xx}^0 = n^2 - 2k + \frac{2\mu}{\left(x + \mu - 1\right)^3} \qquad \text{---- (7)}$$

$$\Omega_{yy}^0 = n^2 - 2k + \frac{\mu}{\left(x + \mu - 1\right)^3} \qquad \text{---- (8)}$$

Here superscript '0' denotes the corresponding value at equilibrium point

Equations (6) and (7) may be written as

$$\Omega_{xx}^0 = B + 2A \quad , \qquad \Omega_{yy}^0 = B - A$$

where

$$A = \frac{\mu}{(x + \mu - 1)^3} \,, \qquad B = n^2 - 2k$$

Let the third body of infinitesimal mass be displaced to $(x_0 + \xi, y_0 + \eta, z_0 + \zeta)$ where $\xi, \eta, \zeta$ are very small quantities and $(x_0 \ y_0 \ z_0)$ is libration point. Using equation (5), the variational equations of motion in the linearized form becomes

$$\ddot{\xi} - 2n\dot{\eta} = \Omega_{xx}^0 \xi \quad , \; \ddot{\eta} - 2n\dot{\xi} = \Omega_{yy}^0 \eta \, , \quad \ddot{\zeta} = \Omega_{zz}^0 \zeta \qquad \qquad ---- (9)$$

Again, using equations (6), (7) and (8) we have

$$\ddot{\xi} - 2n\dot{\eta} = (B + 2A)\xi \,, \qquad \ddot{\eta} + 2n\dot{\xi} = (B - A)\eta \,, \qquad \ddot{\zeta} = -(2k + A)\zeta \qquad \qquad ---- (10)$$

Now, we consider $\quad \xi = A'e^{\lambda t}, \eta = B'e^{\lambda t}, \zeta = C'e^{\lambda t} \quad$ where $A', B', C'$ are very small quantities and $\lambda$ is a small parameter.

Using the above values in (10), we have from 3$^{rd}$ equation of (10)

$$\lambda^2 C'e^{\lambda t} = -(2k + A)C'e^{\lambda t} \quad , \quad \lambda^2 = \pm i(2k + A)$$

which shows that the motion parallel to z axis is stable.

Now, from first two equations of (10), we have,
s
$$\{\lambda^2 - (B + 2A)\}A' - 2nB'\lambda = 0$$
$$2nA'\lambda + \{\lambda^2 - (B - A)\}B' = 0$$

This will have a non trivial solution for A' and B' if

$$\begin{vmatrix} \lambda^2 - (B + 2A) & -2n\lambda \\ 2n\lambda & \lambda^2 - (B - A) \end{vmatrix} = 0$$

i.e. $\quad \lambda^4 - (2B + A - 4n^2)\lambda^2 + (B - A)(B + 2A) = 0 \qquad \qquad ---- (11)$

It is the characteristic equation of problem

$$\lambda^2 = \frac{\left(2B + A - 4n^2\right) \pm \sqrt{D}}{2}$$

where D is the discriminate for critical mass $\mu_c$, $D = 0$

or, $\left(2B + A - 4n^2\right)^2 - 4(B - A)(B + 2A) = 0$

on putting the value of B in above equation we have,

$$\left(2B + A - 4n^2\right)^2 = 4\left(B^2 + AB - 2A^2\right)$$

ie. $\left(2n^2 - 4k + A - 4n^2\right)^2 = 4\left\{\left(n^2 - 2k\right)^2 + A\left(n^2 - 2k\right) - 2A^2\right\}$

ie. $\left(2n^2 + 4k - A\right)^2 = 4\left\{\left(n^2 - 2k\right)^2 + A\left(n^2 - 2k\right) - 2A^2\right\}$

ie. $9A^2 - 8n^2 A + 3n^2 k = 0$

or $A = \frac{8n^2}{18} \pm \frac{8\left(n^2 - 9k\right)}{18}$ ---- (12)

If the positive sign of above equation is considered then,

$$A_+ = \frac{8}{18}n^2 + \frac{8}{18}n^2 - 4k$$

$$= \frac{8}{9} + \frac{4}{3}A_1 - 4k$$ ---- (13)

if negative sign is considered

$A_- = 4k$ ---- (14)

Now, we know that

$$A = \frac{\mu}{\left(x + \mu - 1\right)^3}$$

putting the value of $x$ from equation (4), we have,

$$\mu\left(x_0 + \mu - 1\right)^{-3} = A$$

ie. $-\mu\left\{1 - \left(x_0 + \mu\right)\right\}^{-3} = A$

ie. $-\mu\left\{1 - \left(2\mu - \frac{3}{2}\mu A_1 + 4k\mu\right)\right\}^{-3} = A$

ie. $-\mu\left\{1 + \left(6\mu - \frac{9}{2}\mu A_1 + 12k\mu\right)\right\} = A$

ie.    $-\mu\left\{1+\mu\left(6-\dfrac{9}{2}A_1+12k\right)\right\}=A$

$$\mu_c=\dfrac{-1}{\left(12-9A_1+24k\right)}\pm\dfrac{\left[1-A\left(12-9A_1+24k\right)\right]}{\left(12-9A_1+12k\right)}$$    ---- (15)

taking the positive sign in equation (15)

$$\mu_{c+}=\dfrac{-1}{\left(12-9A_1+24k\right)}+\dfrac{\left[1-A\left(12-9A_1+24k\right)\right]}{\left(12-9A_1+24k\right)}$$

Using the positive value of A is $A_+$

$$\mu_{c++}=-\dfrac{1}{\left(12-9A_1+24k\right)}+\dfrac{\left[1-A_+\left(12-9A_1+24k\right)\right]}{\left(12-9A_1+24k\right)}$$

Using the value of $A_+$ (13), we have

$$\mu_{c++}=-\left(\dfrac{8}{9}+\dfrac{4}{3}A_1-4k\right)<0$$    ---- (16)

Similarly using the value of $A_-$ in place of A, we get

$$\mu_{c+-}=-\dfrac{1}{\left(12-9A_1+24k\right)}+\dfrac{\left[1-A_-\left(12-9A_1+24k\right)\right]}{\left(12-9A_1+24k\right)}$$

Using the value of $A_-$ from (14)

$$\mu_{c+-}=-4k<0$$    ---- (17)

Again taking negative sign in (16)

$$\mu_{c+-}=-\dfrac{1}{\left(12-9A_1+24k\right)}-\dfrac{\left[1-A\left(12-9A_1+24k\right)\right]}{\left(12-9A_1+24k\right)}$$

Using the value of $A_+$ in place of A in the above equation

$$\mu_{c-+}=-\dfrac{1}{\left(12-9A_1+24k\right)}-\dfrac{\left[1-A_+\left(12-9A_1+24k\right)\right]}{\left(12-9A_1+24k\right)}$$

$$\mu_{c-+}=-\dfrac{2}{\left(12-9A_1+24k\right)}+\left(\dfrac{8}{9}+\dfrac{4}{3}A_1-4k\right)$$    ---- (18)

Again using the value of $A_-$ in the equation in place of A

$$\mu_{c--} = -\frac{1}{(12-9A_1+24k)} - \left[ \frac{1-A_-(12-9A_1+24k)}{(12-9A_1+24k)} \right]$$

$$\mu_{c--} = -\frac{2}{(12-9A_1+24k)} + 4k \qquad \text{---- (19)}$$

which is a very small quantity and liable to be neglected.

Thus we get the admissible value of critical mass $\mu_{c++}$ as deduced in equation (16s)

$$\mu_c = -\left( \frac{8}{9} + \frac{4}{3}A_1 - 4k \right)$$

Putting $A_1 = 0$ & $k = 0$ we get $\mu_c = -\frac{8}{9}$ which agrees with the result of Robe (1977) the (-) sign appears because the value of $\mu$ lies on the negative side of origin in the $x-$ axis as investigated in the orientation of the configuration.

Thus $\mu_c = -\mu_0 - \frac{4}{3}A_1 + 4k$ where $\mu_0 = \frac{8}{9}$ and at $A_1 = 0$ & $k = 0$

we have, $\mu_c = -\mu_0$

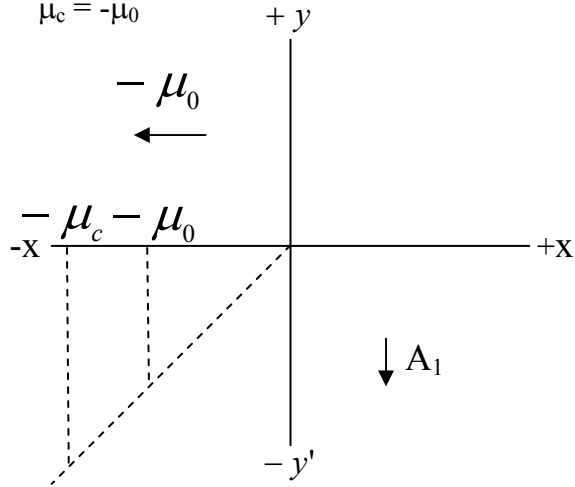

If $-\mu_c < -\mu \leq 0$ _or_ $\mu_c > \mu \geq 0$ or $0 \leq \mu < \mu_c$

we have a range of stability

And if $-\frac{8}{9} \leq -\mu < \mu_c$ or $\frac{8}{9} \geq \mu \geq \mu_c$ or $\mu_c \leq \mu \leq \frac{8}{9}$

we have a range of instability.

**References**

[1.]. **Hallan, P.P. and Neelam Rana (2001)**: The existence and stability of equilibrium points in the Robe's restricted three body problem. Celest. Mech. & Astr. Dyn., $145 - 155$.

[2.] **Kumar, V. and Choudhry, R.K. (1987)**: On the stability of the triangular libration points for the photogravitational circular restricted problem of three bodies when both the attracting bodies are radiating as well. Indian J. Pure Appl. Math. **18(11**), 1023-1039.

[3.] **Liapunov, A.M. (1956)**: A geneal problem of stability of motion Acad. Sc. USSR.

[4.] **Markellos, V.V., Perdios, E. and Labropoulou, P. (1992)**: Linear stability of the triangular equilirbium points in the photogravitational elliptic restricted problem I. Astrophy. and Space Sci. **19**4, 2, 207-213.

[5.] **Mc Cuskey, S.W. (1963)**: Introduction to Celestial Mechanics, Addison-Wesley Publishing Company, Inc., New York.

[6.] **Mishra, P. and Ishwar, B. (1995)**: Second order normalization in the generalized restricted problem of three bodies, small primary being an oblate spheroid. Bull. Astr. Soc. (India) **23**, 4

[7.] **Robe H.A.G. (1978)**: A new kind of three – body problem. Celestial Mechanics, 16, $343 - 351.$ ss

[8.] **Sahoo, S.K. and Ishwar, B. (2000):** Stability of collinear equilibrium points in the generalised photogravitational elliptic restricted three body problem. Bull. Astr. Soc. (India) **2**8, 576-586.

[9.] **Sharma, R.K. (1987):** The linear stability of libration points of the photogravitational restricted three body problem when the smaller primary is an oblate spheroid. Astrophy. and Space Sci. **13**5, 271-281.

[10.] **Shrivastaba, A. K. and Garain, D. (1991)**: Effect of perturbation on the location of libration point in the Robe restricted problem of three bodies (1991) Celest. Mech. & Dyn. Astr. **51,** $67 - 73.$

[11.] **Singh, J. and Ishwar, B. (1999)**: sStability of triangular points in the generalized photogravitational restricted three body problem. Bull. Astr. Soc. (India) **2**7, 415-424.